%%%%%%%%%%%%%%%%%%%%%%%%%%%%%
%Hyperplane Arrangement Cohomology and Monomials in
%            the Exterior Algebra
%
% David Eisenbud, Sorin Popescu, Sergey Yuzvinsky
%
% 12/30/2000
%
%%%% PlainTex. requires diagrams.tex
% begin.tex
%=========================================================================%
% load begin.tex only once, but keep count to match \bye commands
%=========================================================================%

\ifx\begin\undefined\else\global\advance\srcdepth by
1\expandafter\endinput\fi

\def\begin{}
\newcount\srcdepth
\srcdepth=1

\outer\def\bye{\global\advance\srcdepth by -1
  \ifnum\srcdepth=0
    \def\endcmd{\vfill\eject\nopagenumbers\par\vfill\supereject\end}
  \else\def\endcmd{}\fi
  \endcmd
}

%=========================================================================%
% initialize TeX
%=========================================================================%

\magnification=\magstephalf
\baselineskip=13pt
\hsize = 5.5truein
\hoffset = 0.5truein
\vsize = 8.5truein
\voffset = 0.2truein
\emergencystretch = 0.05\hsize

\newif\ifblackboardbold

% comment out the following line if AMS msbm fonts aren't available
\blackboardboldtrue

%=========================================================================%
% select fonts
%=========================================================================%

\font\titlefont=cmbx12 scaled\magstephalf
\font\sectionfont=cmbx12

% Establish AMS blackboard bold fonts without using amssym.def, amssym.tex

\newfam\bboldfam
\ifblackboardbold
\font\tenbbold=msbm10
\font\sevenbbold=msbm7
\font\fivebbold=msbm5
\textfont\bboldfam=\tenbbold
\scriptfont\bboldfam=\sevenbbold
\scriptscriptfont\bboldfam=\fivebbold
\def\bbold{\fam\bboldfam\tenbbold}
\else
\def\bbold{\bf}
\fi

%=========================================================================%
% font size-changing command ("A Beginner's Book of TeX" p35, p275)
%=========================================================================%

\font\Arm=cmr8
\font\Ai=cmmi8
\font\Asy=cmsy8
\font\Abf=cmbx8
\font\Brm=cmr6
\font\Bi=cmmi6
\font\Bsy=cmsy6
\font\Bbf=cmbx6
\font\Crm=cmr5
\font\Ci=cmmi5
\font\Csy=cmsy5
\font\Cbf=cmbx5

\ifblackboardbold
\font\Abbold=msbm10 at 8pt
\font\Bbbold=msbm7 at 6pt
\font\Cbbold=msbm5
\fi

\def\smallmath{%
\textfont0=\Arm \scriptfont0=\Brm \scriptscriptfont0=\Crm
\textfont1=\Ai \scriptfont1=\Bi \scriptscriptfont1=\Ci
\textfont2=\Asy \scriptfont2=\Bsy \scriptscriptfont2=\Csy
\textfont\bffam=\Abf \scriptfont\bffam=\Bbf \scriptscriptfont\bffam=\Cbf
\def\rm{\fam0\Arm}\def\mit{\fam1}\def\oldstyle{\fam1\Ai}%
\def\bf{\fam\bffam\Abf}%
\ifblackboardbold
\textfont\bboldfam=\Abbold
\scriptfont\bboldfam=\Bbbold
\scriptscriptfont\bboldfam=\Cbbold
\def\bbold{\fam\bboldfam\Abbold}%
\fi
}

%=========================================================================%
% single-pass symbolic theorem labeling
%=========================================================================%

% Because this is a single-pass mechanism with no .aux file, forward
% references need to be declared in advance:

%   \forward{thm:main}{Theorem}{1.1}

% This is also the mechanism for "timely" declaration of labels, which
% will usually be buried within the corresponding theorem macros.
% A warning is issued if a label redeclaration is inconsistent, allowing
% forward references to be manually fixed.

%   \ref{thm:main} produces "Theorem~1.1"
%   \refs{thm:main} produces "Theorems~1.1"
%   \refn{thm:main} produces "1.1"

% Some TeX adapted from "The Advanced TeXbook" by David Salomon, chapter 9.

% Implementers: The code for \forward is subtle. Its second argument must
% be provided literally, e.g. "Theorem" rather that "\capitalize{theorem}".
% Its third argument must either be literal or a macro that expands
% directly to a literal, e.g. "\edef\numtoks{\number\proccount}".
% This use of \edef cannot be replaced by \def, which defers expansion.
% Failure to follow these rules will cause spurious warnings that forward
% references are inconsistent, when they are in fact consistent after
% expansion. Note the "Towers of Palo Alto" recreational math problem
% involving the iterated use of \expandafter to expand the first argument
% to \forwardsub before calling it.

\newlinechar=`@
\def\forwardmsg#1#2#3{\immediate\write16{@*!*!*!* forward reference should
be: @\noexpand\forward{#1}{#2}{#3}@}}
\def\nodefmsg#1{\immediate\write16{@*!*!*!* #1 is an undefined reference@}}

\def\forwardsub#1#2{\def\newref{{#2}{#1}}}

\def\forward#1#2#3{%
\expandafter\expandafter\expandafter\forwardsub\expandafter{#3}{#2}
\expandafter\ifx\csname#1\endcsname\relax\else%
\expandafter\ifx\csname#1\endcsname\newref\else%
\forwardmsg{#1}{#2}{#3}\fi\fi%
\expandafter\let\csname#1\endcsname\newref}

\def\firstarg#1{\expandafter\argone #1}\def\argone#1#2{#1}
\def\secondarg#1{\expandafter\argtwo #1}\def\argtwo#1#2{#2}

\def\ref#1{\expandafter\ifx\csname#1\endcsname\relax
  {\nodefmsg{#1}\bf`#1'}\else
  \expandafter\firstarg\csname#1\endcsname
  ~\expandafter\secondarg\csname#1\endcsname\fi}

\def\refs#1{\expandafter\ifx\csname#1\endcsname\relax
  {\nodefmsg{#1}\bf`#1'}\else
  \expandafter\firstarg\csname #1\endcsname
  s~\expandafter\secondarg\csname#1\endcsname\fi}

\def\refn#1{\expandafter\ifx\csname#1\endcsname\relax
  {\nodefmsg{#1}\bf`#1'}\else
  \expandafter\secondarg\csname #1\endcsname\fi}

%=========================================================================%
% widow control
%=========================================================================%

% usage:
% \widow{.2} % start new page if <.2 page left

\def\widow#1{\vskip 0pt plus#1\vsize\goodbreak\vskip 0pt plus-#1\vsize}

%=========================================================================%
% sections and theorems
%=========================================================================%

% use \showlabels or \showlabelsabove to display section and theorem labels

\def\marginlabel#1{}

\def\showlabelsabove{
\font\labelfont=cmss10 at 6pt
\def\marginlabel##1{\rlap{\smash{\raise 10pt\hbox{\labelfont##1}}}}
}

\newcount\seccount
\newcount\proccount
\seccount=0
\proccount=0

\def\stdskip{\vskip 9pt plus3pt minus 3pt}
\def\stdbreak{\par\removelastskip\penalty-100\stdskip}

\def\proof{\stdbreak\noindent{\sl Proof. }}

\def\qed{\vrule height 1.2ex width .9ex depth .1ex}

\def\Box{
  \ifmmode\eqno\qed
  \else\ifvmode\removelastskip\line{\hfil\qed}
  \else\unskip\quad\hskip-\hsize
    \hbox{}\hskip\hsize minus 1em\qed\par
  \fi\stdbreak\fi}

\def\references{
  \removelastskip
  \widow{.05}
  \vskip 24pt plus 6pt minus 6 pt
  \leftline{\sectionfont References}
  \nobreak\stdskip\noindent}

\def\ifempty#1#2\endB{\ifx#1\endA}
\def\makeref#1#2#3{\ifempty#1\endA\endB\else\forward{#1}{#2}{#3}\fi}

\outer\def\section#1 #2\par{
  \removelastskip
  \global\advance\seccount by 1
  \global\proccount=0\relax
		\edef\numtoks{\number\seccount}
  \makeref{#1}{Section}{\numtoks}
  \widow{.05}
  \vskip 24pt plus 6pt minus 6 pt
  \message{#2}
  \leftline{\marginlabel{#1}\sectionfont\numtoks\quad #2}
  \nobreak\stdskip}

\def\proclamation#1#2{
  \outer\expandafter\def\csname#1\endcsname##1 ##2\par{
  \stdbreak
  \advance\proccount by 1
  \edef\numtoks{\number\seccount.\number\proccount}
  \makeref{##1}{#2}{\numtoks}
  \noindent{\marginlabel{##1}\bf #2 \numtoks\enspace}
  {\sl##2\par}
  \stdbreak}}

\def\othernumbered#1#2{
  \outer\expandafter\def\csname#1\endcsname##1{
  \stdbreak
  \advance\proccount by 1
  \edef\numtoks{\number\seccount.\number\proccount}
  \makeref{##1}{#2}{\numtoks}
  \noindent{\marginlabel{##1}\bf #2 \numtoks\enspace}}}

\proclamation{definition}{Definition}
\proclamation{lemma}{Lemma}
\proclamation{proposition}{Proposition}
\proclamation{theorem}{Theorem}
\proclamation{corollary}{Corollary}
\proclamation{conjecture}{Conjecture}

\othernumbered{example}{Example}
\othernumbered{remark}{Remark}
\othernumbered{construction}{Construction}
\othernumbered{problem}{Problem}
%=========================================================================%
% enable postscript illustrations using epsf.tex
%=========================================================================%

% Usage:
% \draw{70}{fig}{} % draw fig.eps at 70% scale
% \draw{999}{fig}{} % draw fig.eps scaled to width of page

% Optional third argument can be multiple calls to \figtext; see below.
% More generally, the third argument is read in vertical mode, with the
% reference point at the lower left corner of the eps picture, whose
% dimensions are contained in the dimen registers \drawx and \drawy.
% This enables using TeX to generate the text that goes with the picture.
% To request that the picture be widened to respect the added text, 
% examine and modify the dimen registers \ngap, \egap, \sgap, \wgap.
% This is done automatically by the \figtext macro.

% These macros rely on "epsf.tex" which is the lowest level interface
% available for including encapsulated Postscript files in TeX documents.
% Rather that manually reading the .eps file to compute the nominal size,
% the \epsfbox macro is called twice, and two of its internal registers
% are examined after the first call. A major change to epsf.tex (unlikely)
% will require changes here. 

%\input epsf

\newcount\figcount
\figcount=0
\newbox\drawing
\newcount\drawbp
\newdimen\drawx
\newdimen\drawy
\newdimen\ngap
\newdimen\sgap
\newdimen\wgap
\newdimen\egap

\def\drawbox#1#2#3{\vbox{
  \setbox\drawing=\vbox{\offinterlineskip\epsfbox{#2.eps}\kern 0pt}
  \drawbp=\epsfurx
  \advance\drawbp by-\epsfllx\relax
  \multiply\drawbp by #1
  \divide\drawbp by 100
  \drawx=\drawbp truebp
  \ifdim\drawx>\hsize\drawx=\hsize\fi
  \epsfxsize=\drawx
  \setbox\drawing=\vbox{\offinterlineskip\epsfbox{#2.eps}\kern 0pt}
  \drawx=\wd\drawing
  \drawy=\ht\drawing
  \ngap=0pt \sgap=0pt \wgap=0pt \egap=0pt 
  \setbox0=\vbox{\offinterlineskip
    \box\drawing \ifgridlines\drawgrid\drawx\drawy\fi #3}
  \kern\ngap\hbox{\kern\wgap\box0\kern\egap}\kern\sgap}}

\def\draw#1#2#3{
  \setbox\drawing=\drawbox{#1}{#2}{#3}
  \advance\figcount by 1
  \goodbreak
  \midinsert
  \centerline{\ifgridlines\boxgrid\drawing\fi\box\drawing}
  \smallskip
  \vbox{\offinterlineskip
    \centerline{Figure~\number\figcount}
    \smash{\marginlabel{#2}}}
  \endinsert}

\def\nextfigtoks{%
  \advance\figcount by 1%
  \edef\numtoks{\number\figcount}%
  \advance\figcount by -1}

\newif\ifgridlines
\newbox\figtbox
\newbox\figgbox
\newdimen\figtx
\newdimen\figty

\newdimen\bwd
\bwd=2sp % 2sp (1/32768") is smallest visible width for Textures

\def\hline#1{\vbox{\smash{\hbox to #1{\leaders\hrule height \bwd\hfil}}}}

\def\vline#1{\hbox to 0pt{%
  \hss\vbox to #1{\leaders\vrule width \bwd\vfil}\hss}}

\def\clap#1{\hbox to 0pt{\hss#1\hss}}
\def\vclap#1{\vbox to 0pt{\offinterlineskip\vss#1\vss}}

\def\hstutter#1#2{\hbox{%
  \setbox0=\hbox{#1}%
  \hbox to #2\wd0{\leaders\box0\hfil}}}

\def\vstutter#1#2{\vbox{
  \setbox0=\vbox{\offinterlineskip #1}
  \dp0=0pt
  \vbox to #2\ht0{\leaders\box0\vfil}}}

\def\crosshairs#1#2{
  \dimen1=.002\drawx
  \dimen2=.002\drawy
  \ifdim\dimen1<\dimen2\dimen3\dimen1\else\dimen3\dimen2\fi
  \setbox1=\vclap{\vline{2\dimen3}}
  \setbox2=\clap{\hline{2\dimen3}}
  \setbox3=\hstutter{\kern\dimen1\box1}{4}
  \setbox4=\vstutter{\kern\dimen2\box2}{4}
  \setbox1=\vclap{\vline{4\dimen3}}
  \setbox2=\clap{\hline{4\dimen3}}
  \setbox5=\clap{\copy1\hstutter{\box3\kern\dimen1\box1}{6}}
  \setbox6=\vclap{\copy2\vstutter{\box4\kern\dimen2\box2}{6}}
  \setbox1=\vbox{\offinterlineskip\box5\box6}
  \smash{\vbox to #2{\hbox to #1{\hss\box1}\vss}}}

\def\boxgrid#1{\rlap{\vbox{\offinterlineskip
  \setbox0=\hline{\wd#1}
  \setbox1=\vline{\ht#1}
  \smash{\vbox to \ht#1{\offinterlineskip\copy0\vfil\box0}}
  \smash{\vbox{\hbox to \wd#1{\copy1\hfil\box1}}}}}}

\def\drawgrid#1#2{\vbox{\offinterlineskip
  \dimen0=\drawx
  \dimen1=\drawy
  \divide\dimen0 by 10
  \divide\dimen1 by 10
  \setbox0=\hline\drawx
  \setbox1=\vline\drawy
  \smash{\vbox{\offinterlineskip
    \copy0\vstutter{\kern\dimen1\box0}{10}}}
  \smash{\hbox{\copy1\hstutter{\kern\dimen0\box1}{10}}}}}

\def\figtext#1#2#3#4#5{
  \setbox\figtbox=\hbox{#5}
  \dp\figtbox=0pt
  \figtx=-#3\wd\figtbox \figty=-#4\ht\figtbox
  \advance\figtx by #1\drawx \advance\figty by #2\drawy
  \dimen0=\figtx \advance\dimen0 by\wd\figtbox \advance\dimen0 by-\drawx
  \ifdim\dimen0>\egap\global\egap=\dimen0\fi
  \dimen0=\figty \advance\dimen0 by\ht\figtbox \advance\dimen0 by-\drawy
  \ifdim\dimen0>\ngap\global\ngap=\dimen0\fi
  \dimen0=-\figtx
  \ifdim\dimen0>\wgap\global\wgap=\dimen0\fi
  \dimen0=-\figty
  \ifdim\dimen0>\sgap\global\sgap=\dimen0\fi
  \smash{\rlap{\vbox{\offinterlineskip
    \hbox{\hbox to \figtx{}\ifgridlines\boxgrid\figtbox\fi\box\figtbox}
    \vbox to \figty{}
    \ifgridlines\crosshairs{#1\drawx}{#2\drawy}\fi
    \kern 0pt}}}}

% macros to add space to text on specified sides

\def\hpad#1#2#3{\hbox{\kern #1\hbox{#3}\kern #2}}
\def\vpad#1#2#3{\setbox0=\hbox{#3}\dp0=0pt\vbox{\kern #1\box0\kern #2}}

% macro to give one text string the apparent height of another

% macro to center one text string over another

\def\stack#1#2#3{\vbox{\offinterlineskip
  \setbox2=\hbox{#2}
  \setbox3=\hbox{#3}
  \dimen0=\ifdim\wd2>\wd3\wd2\else\wd3\fi
  \hbox to \dimen0{\hss\box2\hss}
  \kern #1
  \hbox to \dimen0{\hss\box3\hss}}}

% macros to hide size of trailing exponents

\def\hexp#1{%
  \setbox0=\hbox{${}^{#1}$}%
  \hbox to .5\wd0{\box0\hss}}

%=========================================================================%
% macros for matrices and arrows
%=========================================================================%

% typical usage:
%   \rightarrowmat{2pt}{4pt}{d & bd \cr \!-c & 0 \cr 0 & -ac \cr}

\def\bmatrix#1#2{{\smallmath\left[\vcenter{\halign
  {&\kern#1\hfil$##\mathstrut$\kern#1\cr#2}}\right]}}

\def\rightarrowmat#1#2#3{
  \setbox1=\hbox{\kern#2$\bmatrix{#1}{#3}$\kern#2}
  \,\vbox{\offinterlineskip\hbox to\wd1{\hfil\copy1\hfil}
    \kern 3pt\hbox to\wd1{\rightarrowfill}}\,}

\def\leftarrowmat#1#2#3{
  \setbox1=\hbox{\kern#2$\bmatrix{#1}{#3}$\kern#2}
  \,\vbox{\offinterlineskip\hbox to\wd1{\hfil\copy1\hfil}
    \kern 3pt\hbox to\wd1{\leftarrowfill}}\,}

\def\rightarrowbox#1#2{
  \setbox1=\hbox{\kern#1\hbox{\smallmath #2}\kern#1}
  \,\vbox{\offinterlineskip\hbox to\wd1{\hfil\copy1\hfil}
    \kern 3pt\hbox to\wd1{\rightarrowfill}}\,}

\def\leftarrowbox#1#2{
  \setbox1=\hbox{\kern#1\hbox{\smallmath #2}\kern#1}
  \,\vbox{\offinterlineskip\hbox to\wd1{\hfil\copy1\hfil}
    \kern 3pt\hbox to\wd1{\leftarrowfill}}\,}

%=========================================================================%
% quire macros for preview mode and making booklets
%=========================================================================%

% \legalbooklet{20} makes a booklet from legal paper in landscape
% orientation, where "20" is the page count. To preview, give a negative
% pagecount. Either print using the legal duplex option on a modern laser
% printer, or struggle to simulate this effect manually. Bind using a long
% reach stapler.

% \preview squeezes two pages side by side in landscape orientation. It
% is not suitable for printing, but ideal for previewing on a two page
% monitor.

% \twoup squeezes two pages onto letter paper in landscape mode,
% suitable for printing.

% Each of these macros calls the file "quire.tex"

\def\bookletdims{
  \hsize=5.25truein
  \vsize=7truein
}

\def\legalbooklet#1{
  \input quire
  \bookletdims
  \htotal=7.0truein
  \vtotal=8.5truein
  % below computed from above
  \hoffset=\htotal
  \advance\hoffset by -\hsize
  \divide\hoffset by 2
  \voffset=\vtotal
  \advance\voffset by -\vsize
  \divide\voffset by 2
  \advance\voffset by -.0625truein
  \shhtotal=2\htotal
  % below doesn't need to change
  \horigin=0.0truein
  \vorigin=0.0truein
  \shstaplewidth=0.01pt
  \shstaplelength=0.66truein
  \shthickness=0pt
  \shoutline=0pt
  \shcrop=0pt
  \shvoffset=-1.0truein
  \ifnum#1>0\quire{#1}\else\qtwopages\fi
}

\def\preview{
  \input quire
  \bookletdims
  \hoffset=0.1truein
  \vtotal=8.5truein
  \shhtotal=14truein
  % below computed from above
  \voffset=\vtotal
  \advance\voffset by -\vsize
  \divide\voffset by 2
  \advance\voffset by -.0625truein
  \htotal=2\hoffset
  \advance\htotal by \hsize
  % below doesn't need to change
  \horigin=0.0truein
  \vorigin=0.0truein
  \shstaplewidth=0.5pt
  \shstaplelength=0.5\vtotal
  \shthickness=0pt
  \shoutline=0pt
  \shcrop=0pt
  \shvoffset=-1.0truein
  \qtwopages
}

\def\twoup{
  \input quire
  \hsize=4.79452truein % 5.25/1.095
  \vsize=7truein
  \vtotal=8.5truein
  \shhtotal=11truein
  % below computed from above
  \hoffset=-2\hsize
  \advance\hoffset by \shhtotal
  \divide\hoffset by 6
  \voffset=\vtotal
  \advance\voffset by -\vsize
  \divide\voffset by 2
  \advance\voffset by -12truept
  \htotal=2\hoffset
  \advance\htotal by \hsize
  % below doesn't need to change
  \horigin=0.0truein
  \vorigin=0.0truein
  \shstaplewidth=0.01pt
  \shstaplelength=0pt
  \shthickness=0pt
  \shoutline=0pt
  \shcrop=0pt
  \shvoffset=-1.0truein
  \qtwopages
}

%=========================================================================%
% timestamp (adapted from eplain.tex)
%=========================================================================%

\newcount\countA
\newcount\countB
\newcount\countC

\def\monthname{\begingroup
  \ifcase\number\month
    \or January\or February\or March\or April\or May\or June\or
    July\or August\or September\or October\or November\or December\fi
\endgroup}

\def\dayname{\begingroup
  \countA=\number\day
  \countB=\number\year
  \advance\countA by 0 % adjust after each leap day
  \advance\countA by \ifcase\month\or
    0\or 31\or 59\or 90\or 120\or 151\or
    181\or 212\or 243\or 273\or 304\or 334\fi
  \advance\countB by -1995
  \multiply\countB by 365
  \advance\countA by \countB
  \countB=\countA
  \divide\countB by 7
  \multiply\countB by 7
  \advance\countA by -\countB
  \advance\countA by 1
  \ifcase\countA\or Sunday\or Monday\or Tuesday\or Wednesday\or
    Thursday\or Friday\or Saturday\fi
\endgroup}

\def\timename{\begingroup
   \countA = \time
   \divide\countA by 60
   \countB = \countA
   \countC = \time
   \multiply\countA by 60
   \advance\countC by -\countA
   \ifnum\countC<10\toks1={0}\else\toks1={}\fi
   \ifnum\countB<12 \toks0={\sevenrm AM}
     \else\toks0={\sevenrm PM}\advance\countB by -12\fi
   \relax\ifnum\countB=0\countB=12\fi
   \hbox{\the\countB:\the\toks1 \the\countC \thinspace \the\toks0}
\endgroup}

\def\timestamp{\dayname, \the\day\ \monthname\ \the\year, \timename}

%==========================================================================
% macros (specific to this paper)
%==========================================================================

% surround with $ $ if not already in math mode
\def\enma#1{{\ifmmode#1\else$#1$\fi}}
\def\th{{^{\rm th}}}

\def\mathbb#1{{\bbold #1}}
\def\mathbf#1{{\bf #1}}

% blackboard bold symbols
\def\NN{\enma{\mathbb{N}}}

% bold symbols
\def\aa{\enma{\mathbf{a}}}
\def\bb{\enma{\mathbf{b}}}

\def\ee{\enma{\mathbf{e}}}

\def\uu{\enma{\mathbf{u}}}
\def\kVect{\enma{{k{\mathbf{Vect}}}}}

\def\set#1{\enma{\{#1\}}}
\def\setdef#1#2{\enma{\{\;#1\;\,|\allowbreak
  \;\,#2\;\}}}

\def\sf{\mathop{\rm sf}\nolimits}

\def\Ext{\mathop{\rm Ext}\nolimits}
\def\Tor{\mathop{\rm Tor}\nolimits}
\def\link{\mathop{\rm lk}\nolimits}

\def\supp{\mathop{\rm supp}\nolimits}

\def\A{\enma{\cal A}}
\def\B{\enma{\cal B}}

\def\Sym{\mathop{\rm Sym}\nolimits}

\input diagrams
%\showlabels

%==========================================================================
% macros (specific to this paper)
%==========================================================================

% surround with $ $ if not already in math mode
\def\enma#1{{\ifmmode#1\else$#1$\fi}}
\def\th{{^{\rm th}}}
\def\st{{^{\rm st}}}

\def\mathbb#1{{\bbold #1}}
\def\mathbf#1{{\bf #1}}

% blackboard bold symbols
\def\N{\enma{\mathbb{N}}}
\def\NN{\enma{\mathbb{N}}}
\def\CC{\enma{\mathbb{C}}}
\def\Z{\enma{\mathbb{Z}}}
\def\P{\enma{\mathbb{P}}}

% caligraphic symbols
\def\A{\enma{{\cal A}}}

% bold symbols
\def\aa{\enma{\mathbf{a}}}
\def\bb{\enma{\mathbf{b}}}

\def\ee{\enma{\mathbf{e}}}
\def\FF{\enma{\mathbf{F}}}
\def\LL{\enma{\mathbf{L}}}
\def\RR{\enma{\mathbf{R}}}

\def\uu{\enma{\mathbf{u}}}
\def\kVect{\enma{{k{\mathbf{Vect}}}}}

%caligraphic
\def\F{\enma{\cal{F}}}
\def\O{\enma{\cal{O}}}
\font\abst=cmr9
%
%gothic symbols
% Gothic fonts from AMSTeX 
\font\tengoth=eufm10  \font\fivegoth=eufm5
\font\sevengoth=eufm7
\newfam\gothfam  \scriptscriptfont\gothfam=\fivegoth 
\textfont\gothfam=\tengoth \scriptfont\gothfam=\sevengoth
\def\goth{\fam\gothfam\tengoth}
\def \gm {{\goth m}}

\def\set#1{\enma{\{#1\}}}
\def\setdef#1#2{\enma{\{\;#1\;\,|\allowbreak
  \;\,#2\;\}}}

\def\sf{\mathop{\rm sf}\nolimits}

\def\sing{\mathop{\rm sing}\nolimits}
\def\ann{\mathop{\rm Ann}\nolimits}
\def\Hom{\mathop{\rm Hom}\nolimits}
\def\Ext{\mathop{\rm Ext}\nolimits}

\def\Tor{\mathop{\rm Tor}\nolimits}
\def\tor{\mathop{\rm Tor}\nolimits}
\def\H{{\rm H}}
\def\Tor{\mathop{\rm Tor}\nolimits}
\def\ker{\mathop{\rm ker}\nolimits}
\def\link{\mathop{\rm lk}\nolimits}

\def\supp{\mathop{\rm supp}\nolimits}
\def\initial{\mathop{\rm in}\nolimits}

\def\A{\enma{\cal A}}
\def\B{\enma{\cal B}}

\def\Sym{\mathop{\rm Sym}\nolimits}

%%%%%%%%%%%%%%%%%%%%%%%%%%%%%%%%%%%%%%%%%%%%%
% Forward references
%
\forward{coho of hyp arr}{Section}{1}
\forward{OSrank-var}{Section}{2}
\forward{S-module}{Section}{3}
\forward{Appendix}{Section}{5}
\forward {localsys}{Section}{4}
\forward {coho res}{Theorem}{1.1}
\forward {codim of rank var}{Corollary}{2.3}
\forward{generic example continued}{Example}{3.3}
\forward {characterization of generic}{Corollary}{3.6}

\forward {singular vectors of a standard module}{Theorem}{4.1}
\forward {sf bettis exterior}{Proposition}{5.3}
\forward {betti}{Corollary}{5.7}

%%%%%%%%%%%%%%%%%%%%%%%%%%%%%%%%%%%%%%%%%%%%%

\hbox{}
\bigskip\bigskip
\centerline{\titlefont Hyperplane Arrangement Cohomology}
\smallskip
\centerline{\titlefont  and}
\smallskip
\centerline{\titlefont Monomials in the Exterior Algebra
\footnote{$^{*}$}{\rm Mathematics Subject Classification (MSC 2000) numbers: Primary 15A75,
52C35, 55N45; secondary 55N99, 14Q99.}} 
\bigskip\smallskip
\centerline{David Eisenbud, Sorin Popescu, and Sergey Yuzvinsky
\footnote{$^{**}$}{\rm The first two authors are grateful to the NSF for 
support during the preparation of this work. The authors would like 
to thank the Mathematical Sciences Research Institute in Berkeley 
for its support while part of this paper was being written.}}

\bigskip\bigskip
{\narrower
\noindent{\bf Abstract:} \abst 
We show that if $X$ is the complement of a complex hyperplane
arrangement, then the homology of $X$ has linear free resolution as
a module over the exterior algebra on the first cohomology of $X$.
We study invariants of $X$ that can be deduced from this resolution.
A key ingredient is a result of Aramova, Avramov, and Herzog [2000] on
resolutions of monomial ideals in the exterior algebra. We
give a new conceptual proof of this result.

\bigskip}

\noindent Let $X$ be the complement of a complex hyperplane arrangement $\A$.
In this paper we study the singular homology $\H_*(X)$ as a module
over the exterior algebra $E$ on the first singular cohomology
$V:=\H^1(X)$ always with coefficients in a fixed field $K$.  Our first
main result (\ref{coho of hyp arr}) asserts that $\H_*(X)$ is
generated in a single degree and has a linear free resolution; this
amounts to an infinite sequence of statements about
the multiplication in the Orlik-Solomon algebra
$\H^*(X)$. We also analyze other topological examples from the  
point of view of resolutions over the exterior algebra.

In \ref{OSrank-var} we study an invariant of an $E$-module $N$ called
the {\it singular variety\/}, the algebraic subset of $V$ consisting
of those elements $x$ whose annihilator in $N$ is not equal to $xN$.
The singular variety is the same for $N$ and for $N^*$, and thus for
the homology and cohomology of $X$.  Aramova, Avramov and Herzog
[2000] show that the codimension of the singular variety (called by
them the {\it rank variety\/}) gives the rate of growth of the free
resolution of $N$.  We compute the singular variety of 
$\H^*(X)$ (or $\H_*(X)$): it
is a linear subspace of codimension equal to the number of central arrangements
in an expression of $\A$ as a product of irreducible arrangements.

One way to use the linearity of the resolution of $\H_*(X)$ is through
the Bernstein-Gel'fand-Gel'fand correspondence, which produces a
graded module $F(\A)$ over the symmetric algebra on $W:=V^*=\H_1(X)$
corresponding to $X$, whose associated sheaf $\F(\A)$ is supported on
the singular variety.  In \ref{S-module} we compute this invariant of
$X$ (really an invariant of the intersection poset of the
arrangement).  We compute some homological invariants of $F(\A)$ and
we use its properties to show that the cones over a generic hyperplane
arrangement may be characterized as the arrangements for which the defining
ideal (the Orlik-Solomon ideal) of $\H^*(X)$ also has a linear free
resolution.

If $e\in E_1$, then $e$ corresponds to a local system on the complement
of the hyperplane arrangement $\A$ (see Esnault-Schechtman-Viehweg [1992],
Yuzvinsky [1995], and Libgober-Yuzvinsky [2000], Section 4). 
If we set $A=\H^*(X)$ and $\ann_A e=\setdef{a\in A}{ea=0}$, 
then the homology 
$$
\H(e, A):=(\ann_A e)/eA
$$
of the complex
defined through multiplication by $e$ is the cohomology of $X$ with
supports in that local system, for $e$ subject to a certain mild genericity  
condition (Esnault-Schechtman-Viehweg [1992]).
In \ref{localsys} we show how to compute $\H(e,A)$. In terms of the sheaf
$\F(\A)$ we have
$$
\H(e, A)=\Tor^{\O_{\P,e}}_*(\kappa(e),\F(\A)_e)
$$
where $\P=\P(E_1^*)$ is the projective space of one-dimensional
subspaces of $E_1$, and $\kappa(e)$ is the residue field of the local
ring $\O_{\P,e}$ at the point corresponding to $e$.  (We prove in
\ref{singular vectors of a standard module} a corresponding result
more generally, for arbitrary modules $A$ with linear injective
resolution.)  It follows, for example, that when $e$ is singular on
$A$ the module $\H(e,A)$ has nonzero components in every degree up to
the codimension of $F(\A)$, cf. \ref{singular vectors of a standard
module} $b)$.  This generalizes Theorem 4.1 (i) of Yuzvinsky [1995];
see also Libgober-Yuzvinsky [2000].

A key ingredient in the proof of our main theorem is the theorem of
Aramova, Avramov, and Herzog [2000] (later improved by R\"omer [2001])
relating the resolutions of square-free monomial ideals (and some more
general modules) over symmetric and exterior algebras. This allows us to
apply the results on resolutions and Alexander duality due to Eagon and
Reiner [1998].   The proof given by Aramova, Avramov and Herzog depends
on an intricate computation. In \ref{Appendix} we offer a conceptual
description of the relationship which leads to a transparent proof.

We are pleased to acknowledge the essential role of the computer
algebra system written by Grayson and Stillman [Macaulay2] in the
genesis of this paper: It was only through ``playing" with this
program that we were lead to guess at the main result (\ref{coho res})
and most of the other results were checked for plausibility
before we looked for proofs.

\bigskip
\noindent{\bf Notation:}
Throughout this paper, $\A$ will denote an essential affine complex
hyperplane arrangement, that is, a set of $n$ affine hyperplanes in
$\CC^\ell$ whose intersection poset has rank $\ell$. We will denote the
complement of the union of the hyperplanes in $\A$ by  $X$. We denote
with $K$  an arbitrary field.

We use notation as in Orlik-Terao [1992].  In particular we write
$A:=A(\A)$ for the Orlik-Solomon algebra of $\A$, isomorphic to the
singular cohomology of $X$ with coefficients in $K$.  
The vector space $E_1=V=\H^1(X)$ has basis $e_1,\ldots,e_n$
corresponding to the hyperplanes of $\A$. Writing
$E$ for the exterior algebra on $E_1$ we have
$A=\H^*(X)=E/I$ where
$I\subset E$ is the {\it Orlik-Solomon ideal\/} generated by the
elements
$$
\partial(e_{i_1}\wedge\cdots\wedge e_{i_t})=
\sum_j (-1)^j
e_{i_1}\wedge\cdots e_{i_{j-1}}\wedge \widehat {e_{i_j}}
                   \wedge e_{i_{j+1}}\cdots \wedge e_{i_t}
$$
for which
$\{ H_{i_1},\dots,H_{i_t}\}$ 
is a minimal linearly dependent set of hyperplanes of $\A$,
and the monomials $e_{i_1}\wedge\cdots\wedge e_{i_t}$ for which
$\{ H_{i_1},\dots,H_{i_t}\}$ 
have empty intersection. For all this see Orlik-Solomon [1980].

We grade $E$ by taking the elements of $V$ to have degree $1$
(this is the opposite convention from that of  Eisenbud and
Schreyer [2000]). The homology module $\H_*(X)$ is dual to $E/I$, and
thus is graded in negative degrees. We denote its component of degree
$-p$ by $\H_p(X)$ ($p\geq 0$).

We will write
$\chi(\A,-)$ for the characteristic polynomial 
of the arrangement $\A$.
For our purposes $\chi$ may be defined by the relation
$$
\chi(\A,t)=t^{\ell}\pi(A,-1/t),
$$ 
where $\pi$ is the Poincar\'e polynomial polynomial of $X$, that is
$$
\pi(t)=\sum_j\dim_K\H^j(X)t^j;
$$ 
see Orlik-Terao [1992, Definition 2.52 and Theorem 3.68].

If $A$ is a skew commutative algebra, we write $A\langle e\rangle$
and $A[t]$  to denote the skew-commutative algebra obtained by
adjoining a variable of degree 1 or 2 respectively; thus $A[t]$
is an ordinary polynomial ring on one commuting variable over
$A$, while if $E$ is the exterior algebra of $V$ then 
$E\langle e\rangle$ is the exterior algebra of $V\oplus Ke$.

\section{coho of hyp arr} The Cohomology of Hyperplane Arrangements

\theorem{coho res} The minimal free resolution of $\H_*(X)$,
regarded as a module over the exterior algebra $E=\wedge(\H^1(X))$
by means of the cap product, has the form
$$
\FF:\qquad  
\dots\to 
E^{\beta_2}(\ell-2)\to 
E^{\beta_1}(\ell-1)\to 
E^{\beta_0}(\ell)\to 
\H_*(X)\to 0.
$$  
The ranks $\beta_i$ may be computed from the formula
$$
\sum_{i=0}^\infty \beta_it^i={(-1)}^{\ell}{\chi(\A,t)\over {(1-t)}^n}
\ \ .
$$

In general we will say that a graded $E$-module $M$ has a {\it linear
resolution\/} if $M$ is generated in a single degree $s$ and has
free resolution with $d^\th$ syzygy
module generated in degree $s+d$; the theorem asserts that $H_*(X)$
has a linear resolution with $s=-\ell$.  

We can interpret 
the statement that a
module of the form $E^{\beta_0}(\ell)$ can map onto $\H_*(X)$ in more
familiar language:

\corollary{socle} An element $c\in \H^*(X)$ is annihilated by
the (cup) product with every element of $\H^1(X)$ if and only if
$c\in \H^\ell(X)$.

\proof
Because $E^{\beta_0}(\ell)$ maps onto $\H_*(X)$, we see that $\H_*(X)$
is generated as an $E$-module by $\H_\ell(X)$. In particular we
recover the well-known fact that $\H^j(X)= (\H_j(X))^*= 0$ for
$j>\ell$, so that every element of $\H^\ell(X)$ is annihilated by
$\H^1(X)$.

Conversely, let $c\in\H^*(X)$ be annihilated by $\H^1(X)$.  The
Orlik-Solomon description shows that $\H^*(X)$ is generated as an
algebra by $\H^1(X)$, so $c$ is annihilated by $\H^+$, the ideal of
elements of positive degree in $\H^*(X)$. In particular, $c \cdot
(\H^+\cdot \H_*(X))=0$.  Because $\H_*(X)$ is generated by
$\H_\ell(X)$ we have $\H^+\cdot \H_*(X)=\sum_{j<\ell}\H_j(X)$. 
Since each $\H_j(X)$ is dual to $\H^j(X)$ by the multiplication pairing,
it follows that $c\in \H^\ell(X)$.\Box

For the proof of \ref{coho res} it is convenient to reduce to
the central case. Recall that an arrangement is {\it central\/} 
if the intersection of its hyperplanes is nonempty.
Given a (not necessarily central) arrangement 
$\A$ of $n$ hyperplanes in $\CC^\ell$, we can projectivize and add the
hyperplane at infinity, to get an arrangement in ${\P}^\ell_\CC$;
the affine cone over this arrangement is a central arrangement 
$\B=c\A$ of $n+1$ hyperplanes in $\CC^{\ell+1}$,
called the {\it cone over} $\A$. Conversely, given a central
arrangement $\B$ of $n+1$ hyperplanes in $\CC^{\ell+1}$, 
and a chosen hyperplane $H$ in it, 
we may form the corresponding arrangement of $n+1$ hyperplanes
in projective $\ell$-space. Removing $H$, we get a (not necessarily
central) arrangement $\A=d\B$ of $n$ hyperplanes in 
$\CC^\ell$, which we call the {\it deconing\/} of $\B$
with respect to $H$. For example, $\A$ is the deconing
of $c\A$ with respect to the ``new'' hyperplane.

The Orlik-Solomon algebras
of $\B$ and $d\B$ are closely related. Topologically,
the complement of
the projective arrangement associated to $\B$ is the same as the
complement of the arrangement associated to any of the deconings of $\B$;
thus the complement of $\B$ is a $\CC^*$-bundle over the 
complement of any of the deconings of $\B$. It follows
that the cohomology algebra of any deconing is
canonically isomorphic to the cohomology algebra of the
complement of $\B$ modulo a degree $1$ form.
The following result gives this identification
algebraically. For
this it is convenient to factor the Orlik-Solomon
relations as products of linear forms:

\proposition{deconing} Suppose $\B=\{H_0,\dots,H_n\}$ 
is a central hyperplane arrangement, 
with Orlik Solomon ideal $I$ in the exterior algebra
$E=K\langle e_0,\dots,e_n\rangle$ whose generators
$e_i$ correspond to the hyperplanes $H_i$. 
Let $E'$ be the subalgebra generated by the differences
$e_i-e_j$. Let $I'\subset E'$ be 
ideal generated by
$$
%\eqalign{
\{ (e_{i_1}-e_{i_2})(e_{i_2}-e_{i_3})
%&
\cdots(e_{i_s-1}-e_{i_s})\mid
%\cr
                %&
H_{i_1},\dots,H_{i_s}\hbox{\ are linearly dependent}\}.
%}
$$
The Orlik-Solomon ideal of $\B$ is $I=I'E$, and
$E/I\cong (E'/I')\langle e_j\rangle$ for any $j$. Furthermore,
if $\A$ is the deconing of $\B$ with respect to
$H_j$, then the Orlik-Solomon algebra of $\A$ is 
$E/(I+(e_j))\cong E'/I'$.

\proof One checks directly that
$
(e_{i_1}-e_{i_2})\cdots(e_{i_s-1}-e_{i_s})= 
{(-1)}^s\partial(e_{i_1}\wedge\cdots\wedge e_{i_s}).
$
It follows that $I=I'E$. The rest
of the statements are consequences.\Box
   
\noindent{\sl Proof of \ref{coho res}.\/} To prove
that the resolution of $\H_*(X)$ is linear, we first
reduce to the central case. By \ref{deconing},
the Orlik-Solomon algebra of $c\A$ is 
$A\langle e_0\rangle=A\otimes_KK\langle e_0\rangle$
as skew-commutative algebras, and it follows that the
free resolution of the homology of the complement of $c\A$
is deduced from that of $\A$ by tensoring over $K$ with
$K\langle e_0\rangle$. In particular, one is linear if and only if the
other one is, and we may assume that $\A$ is central to begin with.

With respect to
the lexicographic order on the monomials of $E$, taking
$e_i<e_j$ if $i<j$, the
initial (largest) terms of the generators for the Orlik-Solomon ideal (as
given in the introduction) are
$$\eqalign{ 
\{  
e_{i_1}\wedge\cdots\wedge e_{i_s} &\mid i_1<\dots <i_s,\ \{H_{i_1},\dots,H_{i_s}\}
\hbox{ is an independent set of}\cr 
&\hbox{hyperplanes, and there exists $i_0<i_1$ such that }\cr 
                     &\{H_{i_0},\dots,H_{i_s}\} 
                        \hbox{ is dependent} 
\}. 
}$$
The subsets
that appear in this expression are exactly the broken circuits
of $c\A$. By Bj\"orner [1982], the monomials that are not
divisible by broken circuits are a basis for $A$. It follows that the
generators of $I$ given in the introduction form a Gr\"obner basis.
Consequently the initial ideal of $I$ is the ideal $I_0$ generated by the
monomials in the display. $I_0$
is the {\it broken circuit ideal\/} of 
the matroid defined by the dependence relations among the hyperplanes
of $\A$. 

{}From the general theory of Gr\"obner bases (as for example
in Eisenbud [1995] where the completely parallel theory is treated for
ideals in a polynomial ring) we see that $I_0$ is a flat degeneration
of $I$. More formally, there is an ideal 
$I_t\subset K[t]\otimes_K E$ 
such that the algebra 
$ K[t]\otimes_K E/I_t$ is free (and thus flat) over $ K[t]$, and
$I_0:=(I+(t))/(t)\subset  K[t]\otimes_K E/(t)=E$
is the initial ideal $\initial(I)$, while for $0\neq a\in  K$ we have
$I_a:=(I+(t-a))/(t-a)\subset  K[t]\otimes_K E/(t-a)=E$
is conjugate to $I$ by a linear automorphism of $E$.

The module structure on $\H_*(X)$ comes
from the identification $\H_*(X)=\Hom_K(\H^*(X), K))$,
so $\H_*(X)$ degenerates flatly to
$M_0=\Hom_K(E/I_0,  K)$. More formally, the module
$M=\Hom_{ K[t]}(E/I_t,  K[t])$ is free (and thus flat)
over $ K[t]$, and has special fiber $M/(t)M\cong M_0$,
whereas for $a\neq 0$ the fiber $M_a:=M/(t-a)M$ is conjugate to 
$\H_*(X)$ by an automorphism of $E$.

The first statement of \ref{coho res} amounts to saying that the $k^\th$ graded
component, $\Tor_j^E(\H_*(X), K)_k,$ vanishes for all $j>0$ and 
$k\neq \ell-j$. The vanishing of any one of these vector spaces
is an open condition in flat families, so it suffices to show that 
$M_0=\Hom_K(E/I_0,  K)$ satisfies the conditions of 
\ref{coho res}. 

The algebra $E$ is Gorenstein (injective as a module over itself) with
socle in degree $n$, so
$M_0=\Hom_K(E/I_0,  K)=\Hom_E(E/I_0, E)(n)$ as $E$-modules. 
On the other hand $\Hom_E(E/I_0, E)$ may be identified with the
annihilator $J_0$ of $I_0$, and we see that it suffices to show
that $J_0$ has free resolution of the form
$$
\FF(-n):\qquad  
\dots\to 
E^{\beta_2}(\ell-2-n)\to 
E^{\beta_1}(\ell-1-n)\to 
E^{\beta_0}(\ell-n)\to 
J_0\to 0.
$$

Since $I_0$ is generated by monomials, so is the ideal $J_0$.
Following Aramova, Avramov, and Herzog [2000] (see also 
\ref{Appendix} below for more details) we let $_SI_0$ and $_SJ_0$
be the ideals of $S= K[e_1,\ldots,e_n]$ generated by the monomials
corresponding to the generators of $I_0$ and $J_0$, respectively, so
that $_SI_0$ and $_SJ_0$ are square-free monomial ideals of $S$.
Aramova, Avramov, and Herzog [2000] show that $J_0$ has a free
resolution as above with $d^\th$ syzygies generated in degree
$d+n-\ell$ if and only if $_SJ_0\subset S$ has a resolution with this
same property; for another proof, see \ref{Appendix}, below.

Any square-free monomial ideal $J$ corresponds to a simplicial
complex $\Delta(J)$. Since $I_0$ and $J_0$ are
annihilators of one another in $E$, the simplicial complex
$\Delta(I_0)$ is the Alexander dual of $\Delta(J_0)$; that is, the
faces of $\Delta(J_0)$ are the complements of the nonfaces of
$\Delta(I_0)$. By Eagon-Reiner [1998], $_SJ_0$ has a (linear) resolution 
as above if and only if $_SI_0$ has codimension $n-\ell$ and $S/(_SI_0)$ is
Cohen-Macaulay or, in combinatorial terms, the simplicial complex
$\Delta(I_0)$ is Cohen-Macaulay of dimension $\ell-1$. (See also the
later papers of Terai [1997], Bayer-Charalambous-Popescu [1999],
Musta\c t\v a [2000], and Yanagawa [2000] for
more sophisticated versions of this result.)

It was observed by Hochster [1972] and Stanley [1975]
that the Cohen-Macaulay property of
a simplicial complex follows from a simpler geometric property called {\it
shellability\/}; see also Stanley [1996, Theorem 2.5],  
Bruns-Herzog [1993, Theorem 5.1.13].  It is known that the
simplicial complex corresponding to the broken circuits of a matroid
of rank $\ell$ is shellable of dimension $\ell-1$ (Provan [1977];
see Bj\"orner [1992, 7.4.2(ii) and 7.4.3] and his reference
Billera and Provan [1980]),
concluding the proof of the first statement.

In order to prove the second statement we note that, from the
given resolution, 
$$\eqalign{
\pi(\H_*(X), t)&=\sum_i (-1)^i\pi(E^{\beta_i}(\ell-i),t)\cr
                  &=\sum_i (-1)^i\beta_it^{-\ell+i}(1+t)^n.
}$$                  
On the other hand, since homology and cohomology are dual,
$
\pi(\H_*(X), t)=\pi(A,1/t)=(-1)^\ell\chi(A,-t)/t^\ell,
$
whence the desired formula.\Box

In general we do not know how to write the free resolution of the
Orlik-Solomon ideal explicitly; this seems an interesting problem. 

\remark{surfaces} Here are a few other topological examples
treated from the point of view of resolutions over
the exterior algebra: \hfill\smallskip
$a)$ Perhaps the most familiar topological spaces with
cohomology generated in degree one are compact orientable
surfaces. If $Y$ is an
orientable compact connected surface of genus
$g>0$, then the homology $\H_*(Y)$ {\it does not\/}
satisfy \ref{coho res}, though it comes close: 
by Poincar\'e duality the homology
$\H_*(Y)$ is isomorphic as a module over $E=\wedge \H^1(Y)$
to $\H^*(Y)$, which
has relations of degree $>1$. 
However, if we write $\H^*(Y)=E/I$ then one can check that
(with respect
to any monomial order on $E$) the initial ideal of $I$
is the square-free
stable ideal consisting of all but the last monomial of degree 
$2$ in $E$.
By Aramova, Herzog, Hibi [1998, Corollary
2.5], the initial ideal, and with it
$I$ itself, has linear resolution.  It follows that the  minimal free
resolution of the homology module has the form
$$
\FF:\qquad  
\dots\to 
E^{\beta_3}(-2)\to 
E^{\beta_2}(-1)\to 
E^{{2g\choose 2}-1}\to 
E(2)\to 
\H_*(Y)\to 0.
$$
\smallskip 
$b)$ A result analogous to
\ref{coho res} holds for the homology module of an essential
arrangement of real subspaces of codimension two in ${\mathbb R}^{2n}$
with even dimensional intersections. In this case the cohomology ring of
the complement has again the shape of an Orlik-Solomon algebra. However,
in contrast with the complex case, it is not determined merely by the
intersection lattice,  but requires the knowledge of extra information
on sign patterns (computed as determinants of linear relations, or as
linking numbers in the sense of knot theory); see Bj\"orner-Ziegler
[1992] and Ziegler [1993] for details.\hfill\smallskip  
$c)$ The complements of codimension two subspace arrangements in 
${\mathbb R}^{4}$ 
are equivalent to the link complements obtained
by intersecting them with the three-sphere $S^3$.
More generally, consider
the case of an arbitrary tame link $L=\cup_{i=1}^n L_i$ in $S^3$,
and let $X$ be the compact manifold with boundary that
is the complement of a tubular neighborhood of $L$.
Alexander duality gives
$\dim
\H^1(X)=n$,
$\dim \H^2(X)=n-1$ and $\H^{\ge 3}(X)=0$.
More explicitly, let $e_i\in\H^1(X)$ be the dual of the meridian
of the $i^\th$ boundary component, and let
$f_{i,j}\in\H^2(X)$ be the Alexander dual of the (relative)
homology class of an arc $\gamma_{i,j}$ connecting
the $i^\th$ and $j^\th$ components of the boundary.
The elements $e_i$ form a basis of $\H^1(X)$ and
(with the conventions $f_{i,i}=0$ and $f_{i,j}=-f_{j,i}$)
the $f_{i,j}$ generate $\H^2(X)$.

A Mayer-Vietoris argument shows
that the cohomology ring of $X$ has a presentation
$$
\H^*(X) = \wedge V'/(e_i\wedge e_j - l_{i,j} f_{i,j},\,
f_{i,j} + f_{j,k} + f_{k,i},\, e_k\wedge f_{i,j},\, f_{i,j}\wedge f_{k,l}),
$$
where $V'=\H^1(X)\oplus \H^2(X)$,
the numbers $i$ and $j$ run from $1$ to $n$,
 and 
$l_{i,j}:=\link(L_i,L_j)$  is the linking number of $L_i$ and
$L_j$. In particular, the cohomology algebra
$\H^*(X)$ depends only on the linking numbers 
(for most of this, see Milnor [1957]).

Let $G$ be the graph whose vertices are
the components $L_i$, $i=1,\ldots, n$, and where
two vertices $L_i$ and $L_j$ are connected by an edge if their
linking number $l_{i,j}$ is non-zero. Assume that $G$
is connected and the ground field has characteristic $0$. 
The given relations then suffice to eliminate all the $f_{i,j}$,
and it
follows that
$\H^*(X)$ is generated in degree $1$ (see also
Massey-Traldi [1986, Theorem 1 and Proposition 4.1],
or Matei-Suciu [2000]).

Under these hypotheses, the cohomology ring behaves very
nicely:

\theorem{links} 
If $G$ is connected and the ground field has characteristic $0$, then
both the homology module $\H_*(X)$ and the presentation ideal $I$ of
the cohomology ring $\H^*(X)$ have linear free resolutions over the
exterior algebra $E=\wedge\H^1(X)$.

\noindent{\sl Proof Sketch.}
With these hypotheses the presentation ideal $I\subset E$ is 
generated by the monomials 
$e_i\wedge e_j$, where $i$ and $j$ are vertices not connected by an 
edge in $G$, together with elements
$\sum_k (1/l_{i_k,i_{k+1}})e_{i_k}\wedge e_{i_{k+1}}$, where the
sums is over a cycle in the graph $G$. In particular, $E/I$ is a
quotient ring of the  (exterior algebra) Stanley-Reisner ring of the
graph $G$, regarded as $1$-dimensional simplicial complex on the vertex
set $\{e_1,\ldots,e_n\}$.
  
Now suppose we have chosen $T$ a spanning tree of the connected graph
$G$,  and a total order on the edges of $G$. Recall that an edge $e\in
G\setminus T$ is called {\it externally active} in $T$ if it is the
largest edge in the unique cycle $C_e$ contained in $T\cup \{e\}$. It is
a standard fact that for each enumeration of the edges of $G$ (say
corresponding to  the choice of a monomial order in $E$) there exists a
spanning tree $T_0$ of $G$ such that every edge of $G$ not in $T_0$ is
externally active in $T_0$ (see Bollob\'as [1998, proof of Theorem 10,
p. 351 and Exercise 8, p. 372]). Since the  cycles $C_e$
form a basis of the cycle space of $G$ (see for
example Bollob\'as [1998, proof of Theorem 9, p. 53]), it follows that  
the ideal $I$ has an initial ideal $I_0$, which is the Stanley-Reisner 
ideal of the chosen spanning tree $T_0$ in $G$. 

The fact that the Stanley-Reisner ideal
$I_0$ has a linear resolution follows from Hochster's
formula for the Betti numbers of a square-free monomial ideal (see
Hochster [1977] or, for an exposition, Stanley [1996]) since any
subcomplex of a tree is a forest, which is acyclic in all positive
homological degrees. The  linearity of the injective resolution of
$\H^*(X)$ follows from the fact that 
$T$ is a Cohen-Macaulay simplicial complex as in the proof of
\ref{coho res}.\Box

\section{OSrank-var} The singular variety of an Orlik-Solomon algebra

An element $x\in V= E_1$ is said to be singular on a module $M$ if the
set of elements of $M$ annihilated by $x$ is not the same as $xM$.
The set $\sing(M)$ of singular elements is an algebraic subset of $V$
called the {\it singular variety\/} $\sing(M)$ of $M$; 
see Aramova, Avramov, and
Herzog [2000] for a discussion (but note that the term 
{\it rank variety\/} is used in place of singular variety.) 
These authors prove, among other
things, that the dimension of $\sing(M)$ is
the {\it complexity\/} of $M$, defined as the exponent of growth of
the Betti numbers of $M$. This complexity plays, for
modules over an exterior algebra, a role analogous to that of the
projective dimension for modules over a polynomial ring. In this
section we will compute the singular variety of the Orlik-Solomon algebra $A$
of an arrangement $\A$. It follows at once from the definition that
the singular variety of a module $M$ is the same as that of $\Hom_K(M,
K)$, so this also gives the singular variety of $\H_*(X)$.

Recall that the product $\A_1\times\A_2$ of
arrangements $\A_i$ in $ \CC^{\ell_i}$ is the arrangement $\A$
in $\CC^{\ell_1+\ell_2}$ consisting of the hyperplanes
$H\times \CC^{\ell_2}$ for $H\in \A_1$ and the hyperplanes
$\CC^{\ell_1}\times H$ for $H\in \A_2$. Any arrangement can
be expressed uniquely as the product of irreducible arrangements.
The following well-known remark shows that to compute the singular variety
of the Orlik-Solomon algebra as a module over $E$, it suffices to
treat the irreducible case:

\proposition{tensor product} The Orlik-Solomon algebra of a product
$\A=\A_1\times\A_2$ of two arrangements is given by
$A(\A)=A(\A_1) \otimes_K A(\A_2)$, the tensor product
in the category of graded skew-commutative $K$-algebras. Thus
$\sing(A(\A))=\sing(A(\A_1))\times \sing(A(\A_2))$.

\proof A minimal dependent set of hyperplanes
in $\A$, or a minimal set with empty intersection, comes from
a similar set either in $\A_1$ or in $\A_2$,
proving the first statement. The second follows because 
$A(\A)_1$ is the direct sum of the corresponding spaces for
$\A_1$ and $\A_2$. A linear form $x=(x_1,x_2)$ is singular for
$A(\A)$ if $x_i$ is singular on $A(\A_i)$ for both $i=1,2$.\Box

The main result of this section is:

\theorem{rank of irr} Let $\A$ be an irreducible 
hyperplane arrangement with Orlik-Solomon algebra
$A$ and elements $e_i\in V:=A_1$ corresponding to the hyperplanes
of $\A$. 
\item{$a)$} If $\A$ is noncentral then the singular variety of
$A$ is $V$.
\item{$b)$} If $\A$ is central then the singular variety of
$A$ is the hyperplane spanned by the elements $e_i-e_j$.

\proof
If the singular variety of the Orlik-Solomon algebra $A$ of
an arrangement $\A$ does not contain an element $e\in V$,
then $A$ is a free module over the subring $K[e]/e^2$. It
follows that the Poincar\'e polynomial $1+t$
of $K[e]/e^2$ divides the Poincar\'e
polynomial
$\pi(A,t)$  of $A$.

On the other hand, 
Crapo [1967] 
(see also Schechtman-Terao-Varchenko [1995, Sect.~2])
shows that if $\B$ is an irreducible central arrangement with
deconing
$\A$, then
$$
{\pi(A(\B),t)/(1+t)}_{\mid_{t=-1}}\not=0.
$$ 
It follows that in this case the singular variety of $A(\A)$  contains
everything of degree $1$. In particular, if $\A$ is
an irreducible noncentral arrangement, we may 
apply this remark to $\B=c\A$. Part $a)$ now follows from
\ref{deconing}.

If now $\B$ is an irreducible central arrangement, then
the formula $A(\B)=A(d\B)[e]$ from  \ref{deconing} implies that
the singular variety of $A(\B)$ is equal to the singular variety
of $E'/I'\cong A(d\B)$; that is, it consists of precisely the elements
of $V'$ as required.\Box

{}From \ref{rank of irr} and \ref{tensor product} we get the general
case:

\corollary{codim of rank var}
The singular variety of the Orlik-Solomon algebra $A$ of any
arrangement $\A$ is a linear space of $V$ given in dual coordinates to
the canonical basis $\{e_H\}_{H\in\A}$ by the system of equations
$$\sum_{H\in\A_j} x_H=0,$$
for each central factor $\A_j$. In particular, the 
codimension of the singular variety of $A$ equals  
the number of central factors in an irreducible 
decomposition of $\A$.\Box

\example{generic example}
A central arrangement $\A$ 
in $\CC^\ell$ is called {\it generic\/}
if no set of $\ell$ or fewer hyperplanes of $\A$ is dependent.
Analogously, a noncentral arrangement is called generic 
if every set of $\ell+1$ or fewer hyperplanes meet transversely
(in particular, they don't meet if the number of hyperplanes
is $\ell+1$). 
In the generic noncentral case it follows immediately
from the definition that the Orlik-Solomon ideal $I$ is the 
$(\ell+1)^\st$ power $\gm^\ell$ of the maximal ideal $\gm$
of $E$. From \ref{deconing} it follows from this that in the 
generic central case the
Orlik-Solomon ideal is the $\ell^\th$ power of the maximal ideal of the
subalgebra $E'$ generated by the differences $e_i-e_j$ of the
generators of $E$.

The homology module $\H_*(X)$
is, as for every arrangement, given by
$$
\H_*(X)=\Hom_K(E/I, K)\cong \Hom_E(E/I,E(n))=(0:_EI)(n),
$$
the $n^\th$ twist of the annihilator of $I$.
If $I=\gm^\ell$, then $(0:_EI)(n)=\gm^{n-\ell+1}$.
An explicit computation of the resolution of this ideal is given in
terms of Schur functors in  Eisenbud and Schreyer [2000,
Corollary 5.3]; in particular the resolution is linear.

\section{S-module} The module $F(\A)$ 

Let $W=V^*=\H_1(X)$ be the dual vector space to $V$, and let $S=\Sym(W)$
be the symmetric algebra of $W$, a polynomial ring over $ K$.

As usually stated, the Bernstein-Gel'fand-Gel'fand 
correspondence (BGG) is an isomorphism between the derived category
of bounded complexes of coherent sheaves on ${\P}(V^*$) and the
derived category of bounded complexes of finitely generated graded
modules over $E=\wedge V$.  But if one examines the proof one can
extract a functor ${\RR}$ from the category of graded modules
over $S$ and the category of linear
free complexes over $E$, and also a functor $\LL$
from the category of graded $E$-modules to the category
of linear free complexes over $S$. These functors are equivalences of
categories; see Eisenbud and Schreyer [2000, Proposition 2.1].

Starting with a graded $E$-module $P$ the corresponding complex $\LL(P)$
over $S$ is
$$
\cdots \rTo S\otimes P_i\rTo S\otimes P_{i+1} \rTo \cdots
$$
with differential $1\otimes p\mapsto\sum x_i\otimes e_ip$, where
$x_i$ and $e_i$ are dual bases of $W$ and $V$. Starting with
a graded $S$-module $M$ the corresponding complex $\RR(M)$ over $E$ is
$$
\cdots \rTo \Hom_K(E, M_i)\rTo \Hom_K(E, M_{i+1})\rTo \cdots ,
$$
with differential defined similarly.

Starting from a hyperplane arrangement $\A$, we
consider the injective resolution of $A$ as an $E$-module.
Recall that since $E$ is
Gorenstein, injective resolutions over $E$ are simply
the duals (with respect to $E$ or to $ K$) 
of free resolutions. Thus the
injective resolution of $A$ is the $ K$-dual of the free resolution of
$\H_*(X)$. By \ref{coho res}, this free resolution, and with it the 
injective resolution of $A$, is linear. 

Thus we may define
$F(\A)$  to be the graded $S$-module that is mapped by
$\RR$ to the injective resolution of $A$ as an $E$-module.
The reason for choosing the injective
resolution over the free resolution in the definition of $F(\A)$
is to make
$F(\A)$ finitely generated. 

The following result, which is Theorem 3.7 of Eisenbud and Schreyer
[2000] 
allows us to derive some basic properties of $F(\A)$:

\theorem{reciprocity}
If $M$ is a
graded $S$-module and $P$ is a finitely generated graded
$E$-module, then $\LL(P)$ is a free resolution of $M$ 
if and only if $\RR(M)$ is an injective resolution of $P$. 
\Box

\corollary{F properties} 
$F(\A)$ is generated over $S$ in degree $\ell$ and
has linear free resolution equal to $\LL(A)$. In particular,
\item{$a)$}
$F(\A)$ has projective dimension $\ell$ and $\Ext^\ell_S(F(\A), S)= K$.
\item{$b)$}
 The support of $F(\A)$ is a linear space whose codimension is the 
number of central arrangements in an irreducible decomposition of $\A$.
\item{$c)$}
The Hilbert function of $F(\A)$ is 
$$
\sum_{i=0}^{\infty}\dim_K(F(\A)_i)t^i=
(-1)^{\ell}{{\chi(\A,t)}\over{(1-t)^n}}
$$

\proof By
\ref{coho res} the injective resolution of $A$ over $E$, which
is dual to the free resolution of $H_*(X)$, is linear. 
By \ref{reciprocity},
$$
\LL(A):\quad 
0\rTo S\otimes_K A_0\rTo\cdots\rTo 
S\otimes_K A_\ell \rTo F(\A)\rTo 0
$$ 
is a (linear) free resolution of $F(\A)$, proving the first
statement and computing the projective dimension.

$a)$: The degree 0 and 1 parts of $A$ coincide with those
of $E$; thus the left-hand terms of the resolution above
are the same as those in $\LL(E)$, the Koszul complex. This
allows us to compute the $\Ext$ in part a).

$b)$: Aramova, Avramov and Herzog [2000] show in general that
the singular variety of an $E$-module $P$ 
is the support of the $S=\Ext_E^*( K, K)$-module
$\Ext_E^*(P, K)$, which is the same (since $E$ is
Gorenstein) as the support of the module $\Ext_E^*( K,P)$.
By  Eisenbud and Schreyer [2000, Proposition 2.3]
this is the module $F(\A)$.

$c)$: Knowing the free resolution of $F(\A)$ allows us to
compute its Hilbert series, just as in 
the proof of \ref{coho res}.\Box

\example{generic example continued}
If $\A$ is a generic noncentral arrangement of $n$ hyperplanes in
$ K^\ell$, then $A$ is $E/\gm^{\ell+1}$, so the
free resolution of $F(\A)$ is a truncation of the Koszul complex,
and $F(\A)$ is isomorphic to the $(v-\ell)^\th$ syzygy module
of the trivial $S$-module $ K$.\Box

We have already seen that if $\A$ is a generic noncentral arrangement
then the Orlik-Solomon ideal of $\A$ is a power of the maximal
ideal of $E$, and thus has a linear free resolution. We will
show that this property characterizes generic arrangements and
their cones. We begin with a general result characterizing
deformations of powers of the maximal ideal:

\theorem{characterization of powers}
Let $I\subset E$
be an ideal in the exterior algebra. Both $I$ and $(E/I)^*$
admit linear free resolutions if and only if $I$ reduces to
a power of the maximal ideal modulo some (respectively any)
maximal $E/I$ regular sequence of linear forms of $E$.

\proof If $f_1,\dots,f_s\in E_1$ is a regular sequence on
$E/I$ then $I$ and $(E/I)^*=\Hom_K(E/I,K)$ are also free
over $K\langle f_1,\dots,f_s\rangle$. The freeness of
$E/I$ over $K\langle f_1,\dots,f_s\rangle$ implies that, 
the image of $I$ in $E/(f_1,\dots,f_s)$ is
isomorphic to $I/(f_1,\dots,f_s)I$, and also that the dual
of $E/(I+(f_1,\dots,f_s))$ is $(E/I)^*\otimes_E E/(f_1,\dots,f_s)$.
Thus the minimal  free
resolutions of $I$ and $(E/I)^*$ are linear if and
only if the minimal resolutions of 
$I/(f_1,\dots,f_s)I$ and $(E/(I+(f_1,\dots,f_s)))^*$
are linear, and it follows from  Eisenbud and 
Schreyer [2000, Section 5] that if the image of
$I$ in $E/(f_1,\dots,f_s)$ is a power of the maximal ideal,
then the minimal free resolutions of $I$ and $(E/I)^*$ are
linear.

To prove the converse, the argument given above
reduces us to showing, in the case where the singular variety
of $E/I$ is $V$, that if the resolutions of $I$ and $(E/I)^*$
are linear, then $I$ is itself a power of the maximal ideal.

Our hypothesis implies in particular that  module $(E/I)^*$ is generated
in a single degree. It follows by Nakayama's Lemma and duality 
that the socle of $E/I$ (the annihilator in $E/I$ of $\gm$)
is generated in a single degree, say degree $s$. Thus 
$I_j=E_j$ for $j>s$,
and it suffices to show that $I_j=0$ for $j\leq s$.

By \ref{reciprocity}, both $\LL(E/I)$ and 
$\LL(I^*)$ are free resolutions; let $F$ be the
module whose resolution is $\LL(E/I)$. 
By Aramova, Avramov, and Herzog [2000] its support is
the singular variety of $E/I$, that is, $V$.

Duality (into $ K$) over the 
exterior algebra gives an exact sequence
$0\to (E/I)^*\to E^*\to I^*\to 0$. Taking duals commutes
with the functor $\LL$ (up to shifts), 
so we get an exact sequence
of complexes
$0\to \LL(E/I)^*\to \LL(E)^*\to \LL(I^*)\to 0$, where
now the duals denote $\Hom_S(-, S)$. The homology
of $\LL(E/I)^*$ at $S\otimes_K ((E/I)_s)^*$
is $\Ext_S^0(F,S)$, which is nonzero because $F$
has support $V$. It follows from the exact sequence that 
$\LL(I^*)$ has nonzero homology at the term
$S\otimes (I_{s+1})^*$. Since $\LL(I^*)$ is a resolution,
this must be the last term of the complex---that is,
$I_j=0$ for $j\leq s$, as required.
\Box

\example{nontrivial defo} The ideals characterized in
\ref{characterization of powers} include  powers of the maximal ideal 
in subalgebras generated by linear forms (this will be the
case for cones over hyperplane arrangements) but also
many that are {\it not\/} of this form. Here is
the simplest concrete example: Let 
$$
I:=(ab+cd, ac, bc)\ \subset \ E:=K\langle a,b,c,d\rangle .
$$
It is easy to check that the three given quadrics form
a Gr\"obner basis with respect to any order with $ab>cd$.
Since $d$ is a regular element on $E$ modulo the initial
ideal $(ab,ac,bc)$, it follows that $d$ is regular on $E/I$.
It is evident that $I$ reduces modulo $d$ to the square of the maximal 
ideal. To see that $I$ is not the square of the maximal ideal
of any exterior subalgebra on 3 variables, note that the quadrics
in 3 variables are all of rank 2, where as $I$ contains an element
of rank 4 (here the rank is defined via the identification between
elements of $E_2$ and skew-symmetric $4\times 4$ matrices.)

\corollary{characterization of generic} 
The Orlik-Solomon
ideal of $\A$ admits a linear free resolution over $E$ if and only if
$\A$ is obtained by successively coning a generic noncentral 
arrangement.

\proof We have already seen that the property holds for generic
noncentral arrangements. If $I$ is the Orlik-Solomon ideal
of $\A$, then the Orlik-Solomon ideal of $c\A$ in $E[e_0]$
is $IE[e_0]=I\otimes_EE[e_0]$, which has free resolution obtained
from that of $I$ by tensoring with $E[e_0]$; in particular,
the linearity is not affected. 

Deconing $\A$ as many times as possible,  it now suffices to show 
that if $\A$ is noncentral and $I$ has a linear resolution then
$\A$ is generic. Since $\A$ is noncentral it can have
no central factors in its irreducible decomposition, and
thus the singular variety of $\A$ is the whole of the vector space
$V$ of linear forms. 

The theorem now follows from a more general result. Recall 
from Aramova, Avramov, and Herzog [2000] that a sequence
of elements $f_1,\dots,f_s\in E_1$ is called a regular sequence on an
$E$-module
$M$ if $M$ is free over $K\langle f_1,\dots,f_s\rangle$, or
equivalently, if the annihilator of $f_i$ in $M/(f_1,\dots,f_s)M$
is $f_iM/(f_1,\dots,f_s)M$ for every $i$. In this case
the minimal free resolution of $M/(f_1,\dots,f_s)M$ over
$E/(f_1,\dots,f_s)$ is obtained by reducing the minimal 
$E$-free resolution of $M$ modulo $(f_1,\dots,f_s)$. The
length of any maximal regular sequence on $M$ is equal
to the codimension of the singular variety of $M$ in $V$.

\remark{horrocks} \ref{characterization of powers} is 
actually equivalent to 
the Theorem of Horrocks that characterizes the bundle
$\Omega_{\P(W)}^i(i)$ as the unique indecomposable sheaf ${\cal F}$ 
such that the only nonzero intermediate
cohomology of any twist of ${\cal F}$ is $\H^i({\cal F})=K$.
To see this one uses the correspondence between powers of
the maximal ideal of $E$ and the twisted exterior powers of
the cotangent sheaf $\Omega_{\P(W)}^i(i)$, as well as the relation
between resolutions over $E$ and cohomology of sheaves on
$\P(W)$, all explained in  Eisenbud and Schreyer [2000].

\section{localsys} How singular is a singular element?

We keep the notation $W=V^*$ and $S=\Sym(W)$ from the last
section. Let $N$ be a graded $E$-module, such as $\H^*(X)$.  Denote
the homology of $N\rTo^eN\rTo^eN$ by $\H(e,N)$; that is,
$\H(e,N)=(\ann_N e)/eN$, so that $\H(e,N)\neq 0$ if and only if $e$ is
in the singular variety of $N$. In general, it is not easy to say in
which degrees the homology $\H(e,N)$ will occur, but if $N$ has a
linear injective resolution and socle in one degree, say degree $0$, as
is the case of $\H^*(X)(l)$, then the following result shows that
$\H(e,N)\neq 0$ if and only if $e:\ N_1\to N_0$ fails to be
surjective, or equivalently if and only if $e$ annihilates a degree $0$
element of $N^*$.

\theorem{singular vectors of a standard module} 
Let $N$ be a finitely generated $E$-module with socle in degree $0$
having a linear injective resolution, and let $F$ be the 
 $S$-module such that the resolution of $N$ is $\RR(F)$. For
any $e\in V$ we have
\item{$a)$} For each $i$,
$$
\H(e,N)_i=\Tor^{\O_{\P,e}}_i(\kappa(e),\F_e),
$$ where $\kappa(e)$ denotes the residue class field of $\P=\P(W)$ at
$e$, and $\F_e$ denotes the stalk at $e$ of the sheaf $\F$ on $\P$
corresponding to the module $F$.  
\item{$b)$} $\H(e,N)_i\neq 0$ if and
only if $i$ is between $0$ and the projective dimension of $\F_e$ as
an $\O_{\P,e}$-module.  In particular, $\H(e,N)\neq 0$ if and only if
$\H(e,N)_0\neq 0$.  
\item{$c)$} Suppose $e\in \sing(N)=X$ is a generic
point of the component of $X$ in which it lies.  If the codimension of
$X\subset \P$ at $e$ is $c$, then $\H(e,N)_i\neq 0$ precisely for
$0\leq i\leq c$.

\problem{} Characterize, in terms of the arrangement,
the subset of $\sing(\H^*(X))$ consisting
of those $e$ for which the projective dimension of $\F(\A)_e$ is
greater than the codimension of the largest component of 
$\sing(\H^*(X))$ on which $e$ lies. See Yuzvinsky [1995], Falk [1997],
Cohen-Suciu [1999] and Libgober-Yuzvinsky [2000] for related computations.

\proof  Everything follows from the formula given in $a)$.
To prove this formula, recall that if $t\in W$ is any
element outside $e^\perp$, then $\F_e$ may be represented as the
degree $0$ part of the localization, at $(e^\perp)$, 
of the module $F[t^{-1}]$. Thus if we write $R=S[t^{-1}]_0$, then
the Tor we are concerned with may be represented
as $\tor^R_i(R/(e^\perp)_0, F[t^{-1}]_0)$. Since localizing
and taking the degree 0 part are both exact functors, we
may rewrite this as $\tor^S_i(S/(e^\perp)[t^{-1}], F[t^{-1}])_0$.
By \ref{reciprocity}, $\LL(N)$ is a free resolution of $F$,
and thus $\tor^S_i(S/(e^\perp)[t^{-1}], F[t^{-1}])_0$ is
the $i^\th$ homology of $\LL(N)\otimes_S R_e$, where 
$R_e= K[t,t^{-1}]$ has
$S$-module structure coming from the
map $S\to K[t,t^{-1}]$ derived from $e:W\to K=Kt$. 
But $\LL(N)\otimes_S R_e$ is nothing but the complex
$$
\cdots\rTo^{t\otimes e}
R_e\otimes_K N_{i-1}\rTo^{t\otimes e}
R_e\otimes_K N_{i}\rTo^{t\otimes e}
R_e\otimes_K N_{i+1}\rTo^{t\otimes e}\cdots
$$ 
whose degree 0 part is
$$
\cdots\rTo^{t\otimes e}
Kt^{-i-1}\otimes_K N_{i-1}\rTo^{t\otimes e}
Kt^{-i}\otimes_K N_{i}\rTo^{t\otimes e}
Kt^{-i+1}\otimes_K N_{i+1}\rTo^{t\otimes e}\cdots,
$$
a complex computing $\H(e,N)_i$.
\Box

\example{} Each part of \ref{singular vectors of a standard module} 
fails for a module whose resolution is linear for $d$ steps, no
matter how large $d$ is. To see this,
note first that if we have an exact sequence
$$0\to N'\to F\to N\to 0,$$
such that $F$ is free, then the homology under multiplying by $e$
satisfies $\H(e,N') = \H(e,N)(1)$, Thus a module whose
linear syzygy chain ends at a certain point cannot have a
minimal generator annihilated by a linear form! Perhaps the
simplest example is the dual of the $d^\th$ module of syzygies
of $E/(ab)$, where $a,b$ are independent linear forms.
\smallskip

Restating part of this in the case of a hyperplane arrangement,
we get:

\corollary{highesthom}
An element $e\in V$ is singular
for $A$ (or equivalently for $\H_*(X)$) 
if and only if
there is a nonzero element of $\H_\ell(X)$ 
annihilated by $e$.

\proof We have shown that $\H^*(X)$ has socle in degree $\ell$
(note that the signs are opposite to those of Eisenbud and Schreyer [2000])
and linear injective resolution, so $N=\H^*(X)(\ell)$ satisfies
the hypothesis of \ref{singular vectors of a standard module}.
Thus $e$ is singular for $N$ if and only if 
$e: H^{\ell-1}(X)\to H^\ell(X)$ fails to be surjective, or,
dually, $e:H_\ell(X)\to H_{\ell-1}(X)$ fails to be injective.
\Box

\section{Appendix}
Syzygies of Monomial Ideals in the Exterior Algebra

In this section we give a conceptual description and proof of the
correspondence between free resolutions of certain modules over exterior
and symmetric algebras first proved by Aramova, Avramov, and Herzog
[2000] and R\"omer [2001].  The main idea is an isomorphism between
certain subcategories of the categories of modules over these two
algebras. Our approach provides a simple explanation for the shape of the
formula relating the corresponding multigraded Betti numbers 
(\ref{betti}).

Let $V$ be an $n$-dimensional vector space over the field $K$, with
basis $x_1,\ldots,x_n$. We will denote by $S=\Sym(V)$ 
the symmetric algebra over $V$,  which we identify with 
the ring of polynomials over $K$ in the $n$ variables $x_1,\ldots,x_n$,
and by $E=\Lambda(V)$ the exterior algebra of the vector space 
$V$. Both these algebras have a natural $\Z^n$ grading in which each
monomial (product of the $x_i$) generates a homogeneous component.
(Note that in earlier sections we wrote $S=\Sym(W)$, where
$W$ was the dual of $V$. Since we have explicitly chosen
a basis of $V$ we may identify $V$ with $W$.)

We say that a $\Z^n$-graded module $M$ over $E$ or $S$ is
{\it square-free\/} if it admits a free presentation 
$F\rTo G\rTo M\rTo 0$ where each generator of $F$ and $G$ has the 
degree of a square-free monomial. Note that the presentation map
$F\rTo G$ is represented by a matrix whose entries are scalars
times monomials. Examples include the 
Stanley-Reisner rings $S/I$ where $I$ is an ideal generated by
square-free monomials, but also such things as the cokernel of 
the matrix
$$\pmatrix{
x_0 &0 \cr
-x_1 &x_1 \cr
0 & -x_2},
$$
the canonical module of the cone over 3 points in the plane.

There is a 1-1 correspondence between square-free modules over $S$ and
over $E$ obtained by interpreting the presentations as matrices over
$S$ or over $E$; we will write $_SM$ and $_EM$ for the two. 

We can describe the correspondence of resolutions in a simple way
as follows:

Start from a free resolution of a square-free module $_SM$. Replace
each free module in the resolution by a module made from the sum of
the vector spaces of its multihomogeneous elements of square-free
degree. It turns out---this is the main point---that this complex
of vector spaces has the structure both of a complex of
$S$-modules and a complex of $E$-modules. The modules in this
complex are not free, but they have simple and functorial free
resolutions. The free resolutions of the
$E$-modules in the complex fit together to make a double complex,
whose total complex is the minimal free resolution of $_EM$. A
similar procedure allows us to pass in the opposite direction. 

The correspondence described above works,
with appropriate definitions, in a more general setting,
in which $E$ is replaced by one 
of the algebras
$$
R_q:={K\{\,x_1,\ldots , x_n\,\}\over
({\langle 
x_jx_i-qx_ix_j \mid 1\le i<j\le n\rangle
+
\langle (1-q)x_i^2\mid 1\le i\le n\rangle})}
$$ 
where $K\{\,x_1,\ldots , x_n\,\}$ denotes the free
 $K$-algebra on
$x_1,\ldots , x_n$, and $q\ne 0$. We leave the details of this
generalization to the interested reader.

All modules and free resolutions considered will
be assumed $\Z^n$-graded. We identify $\N^n\subset\Z^n$ with the set
of monomials of $S$. By the support of a monomial in either $E$ or $S$,
we will mean the collection of variables present in it.
A square-free monomial (or multidegree) is an element
$\aa\in\{0,1\}^n\subset\NN^n$, so $\supp(\aa)=\setdef
{x_j} {a_j\ne 0}$.  

\medskip\noindent
{\bf Modules With Square-free Presentation}.
The following result is due Bruns and Herzog [1995, Theorem 3.1 a)]:

\proposition{sf bettis} Let $\Gamma$ be any set of monomials of $S$
closed under taking least common multiples. If $M$ is an $S$-module
with generators and relations having degrees in $\Gamma$, then all the
free modules in a minimal free resolution of $M$ have degrees in
$\Gamma$.

\proof 
We give a new proof using Gr\"obner bases, which will
easily extend to give \ref{sf bettis exterior} as well.
Let 
$F\rTo^\phi G\rTo M\rTo 0$
be a $\Z^n$-graded free presentation with degrees of $F$ and $G$ 
in $\Gamma$. 
We may replace $F\rTo^\phi G$ by a map $F'\rTo^{\phi'} G$ 
so that the generators
of $F'$ map to a Gr\"obner basis of $\ker(G\rTo M)$ by using the Buchberger
algorithm; this involves adding free generators whose
degrees are the least common multiples of pairs of generators
already present, and thus still in $\Gamma$.
Schreyer's theorem 
(Eisenbud [1995, Theorem 15.10]) shows that in the symmetric
case the
kernel of $\phi'$ is generated by elements of degrees equal
to the least common multiples of pairs of degrees of generators of $F'$.
%In the exterior case we must also add the

It follows as in Eisenbud [1995, Theorem 20.2],
that the minimal presentation of $M$  has also 
degrees in $\Gamma$, and iterating this process 
we see that the same is true for
the whole syzygy chain.
\Box

If $M$ is a square-free module in the sense above, then
we say that the {\it square-free part\/} of $M$ is the 
module obtained by factoring out all the homogeneous elements
of $M$ with non square-free degrees. Thus for example 
the square-free part of $S$ itself is the factor ring
$R:=S/(x_1^2,\ldots,x_n^2)$.  More generally, if $\aa$ is
any square-free monomial, then
$S(-\aa)$ has square-free part $R/\supp(\aa)(-\aa)$. 

\corollary{sf parts} If $M$ is a square-free module over $S$
then the square-free part of $M$ 
admits a resolution by direct sums of modules of the
form $R/\supp(\aa)(-\aa)$. 
\Box

An analogous result also holds over $E$:

\proposition{sf bettis exterior} If $M$ is a module over $E$ whose
generators and relations have square-free  degrees, 
then the square-free part of
$M$ admits a finite resolution by modules of the form 
$E_\aa:= E/\supp(\aa)(-\aa)$.

\proof Because the generators and relations of $M$ have square-free degrees,
we may write $M$ as the cokernel of a map (always $\Z^n$-homogeneous)
between finite direct sums of modules of the form $E_\aa$, and it thus
suffices to show that the kernel of such a map is generated in
square-free degrees. Using Gr\"obner bases we may reduce as above to
the monomial case. Exactly as in Eisenbud [1995, Lemma 15.1], one
shows that all the relations among monomials are generated by those
determined by the fact that any monomial $\aa$ is annihilated by the
variables in the support of $\aa$, and the two-at-a time relations
coming from the least common multiples (``divided Koszul
relations''). The desired result follows.\Box

\medskip\noindent
{\bf The Common Subcategory}. 
The category of modules
over $E$ and the category of modules over $R=S/(x_1^2,\dots,x_n^2)$ 
have much in common.
We make one such connection precise as follows:

Let $\aa$ and $\bb$ be two monomials in $E$ such that
$\supp(\aa)\subseteq\supp(\bb)$, and let $E_\aa$ and $E_\bb$
be the cyclic $E$-submodules generated by these monomials.
The natural inclusion $E_\bb\subseteq E_\aa\subseteq E$ induces
a functorial commutative diagram
\newarrow{Equals}=====
$$\diagram[midshaft,small]
E/\supp(\aa)(-\aa)&\rTo^\cong& E_\aa&\rIntoA&E\\ 
\uTo^{\cdot \bb\aa^{-1}} &&\uIntoB&&\uEquals\\
E/\supp(\bb)(-\bb)&\rTo^\cong& E_\bb&\rIntoA&E\\ 
\enddiagram
$$
where the horizontal isomorphisms are defined by sending 1
to the distinguished generator, and the upper left
monomorphism is induced by right multiplication in $E$
with $\bb\aa^{-1}$, the signed exterior monomial such that
$(\bb\aa^{-1})\aa=\bb$.

The same commutative diagram holds if we replace
$E$ by $R$, and in fact
identifying square-free monomials in $E$ with the corresponding
monomials in $R$ defines an equivalence of categories. More precisely:

\proposition{equivER} Let $\A$ denote the $K$-additive extension of the
category of $\Z^n$-graded submodules of $R$, with morphisms given
by inclusions, let $\B$ denote the $K$-additive extension of the
category of $\Z^n$-graded submodules of $E$, also with morphisms given
by inclusions, and let $\kVect$ be the category of $K$-vector spaces.
\hfill\break 
The above identification of square-free monomials in $E$ with those of $R$ 
induces an equivalence $\Psi:\A \rTo \B$ of categories 
%$$
%\diagram[midshaft,small]
% \A&&\rTo^\Psi&&\B\\
% &\rdTo&&\ldTo\\ 
% &&\kVect&&\\
% \enddiagram
% $$
whose restriction (via the natural forgetful functors) 
to the underlying $K$-vector spaces is the
identity functor. 
In particular, the functor $\Psi$ preserves acyclic complexes.\Box

Notice that if $\aa$ is a square-free monomial, then the square
free parts of 
$S(-\aa)$ is $R/\supp(\aa)(-\aa)=(\aa)R$, an object of $\A$.
Similarly, the square-free part $E_\aa=E/\supp(\aa)(-\aa)=(\aa)E$
of $E(-\aa)$ is an object of $\B$.

\medskip\noindent
{\bf Resolutions over $S$ and $E$}. 
We let $\A_0$ and $\B_0$ be the
additive subcategories generated by these modules and the inclusion
morphisms $(\aa)R\subset(\bb)R$ and $(\aa)E\subset(\bb)E$ when
$\aa|\bb$ as monomials in $S$.

Certain free complexes over $S$ and $E$ correspond to complexes 
in the categories $\A_0\cong\B_0$. We describe the connection
with $S$ first:

To a given $\Z^n$-graded complex ${\mathbf F}_\bullet$ of free $S$-modules
$$
{\mathbf F}_\bullet: \quad 0  \rTo F_r 
\rTo \ldots \rTo F_1 \rTo  F_0, 
$$ 
with generators in square-free degrees,
we associate a complex $\sf({\mathbf F}_\bullet)$ of
$R$-modules, that we may regard as a complex in $\A_0$. Namely
we define $\sf({\mathbf F}_\bullet)$ as the complex of
square-free degrees of ${\mathbf F}_\bullet$, that is
$$
{\sf({\mathbf F}_\bullet)}_i=
\oplus_{\bb\in\set{0,1}^n}{({\mathbf F}_i)}_\bb,
$$
for all $i$, and where the differentials are induced by the
differentials of the original complex ${\mathbf F}_\bullet$.
It is easy to see that $\sf$ 
defines a functor from the category of $\Z^n$-graded complexes
of free $S$-modules to the category of complexes in $\A_0$.

It follows from \ref{sf bettis}
that ${\mathbf F}_\bullet$ is square-free acyclic (that is it has
no homology in square-free multidegrees) if and only if the complex
$\sf({\mathbf F}_\bullet)$ is acyclic. It is also clear that
${\mathbf F}_\bullet$ is minimal if and only if
$\sf({\mathbf F}_\bullet)$ is minimal.

We have proven:

\proposition{S to A} The functor $\sf$ is an equivalence between the category
of square-free complexes of free $S$-modules and the category of
complexes in $\A_0$. It preserves minimality and acyclicity.

Now we turn to complexes over $E$. For functorial constructions,
we will use the divided power algebra.
If $U$ is a finitely dimensional graded vector space, we write 
$D_l(U)$ for the $l\th$-divided power of $U$.  It is convenient to define
$D_l(U)$ as the dual of the $l\th$-symmetric power of the dual space,
that is $D_l(U) = (\Sym_l(U^*))^*$.  The divided powers $D_l(U)$ have
``diagonal'' maps 
$D_{l+1}(U) \rTo D_l(U)\otimes U$ 
which are the monomorphisms dual to the surjective 
natural multiplication map in the symmetric algebra
$\Sym_l (U^*)\otimes U^* \rTo \Sym_{l+1} (U^*).$ 

We can now go from complexes in the category $\B_0$ to 
free complexes over $E$ using the
Cartan Resolution.

\proposition{equivEE} There exists a functor $\Phi$ from the
category of complexes in $\B_0$ to the category of complexes of free
modules over $E$ whose inverse is obtained by taking
square-free parts. $\Phi$ preserves acyclicity and minimality.
Applied to an acyclic complex in $\B_0$ with homology
$M$, the functor $\Phi$ provides an $E$-free resolution of $M$.

\proof  We first define $\Phi$ on modules in $\B_0$.
It associates to a cyclic module 
$E_\aa\cong E/\supp(\aa)(-\aa)$ 
the (resolution) 
$\Phi(E_\aa):=D(L_\aa)\otimes E$, 
where 
$L_\aa:=\oplus_{x_i\in\supp(\aa)}k(-\ee_i)$ 
is the $\Z^n$-graded subspace of $V$ spanned by 
$\supp(\aa)$, and whose  differentials are induced by 
the diagonals followed by multiplication in $E$. 

More precisely $\Phi(E_\aa)$ is the complex
$$\Phi(E_\aa):\qquad
\ldots\rTo D_2(L_\aa)\otimes E(-\aa)\rTo L_\aa\otimes E(-\aa)\rTo
E(-\aa),$$
which is a minimal free resolution of  the cyclic module $E_\aa$.
We see at
once that $E_\aa$ is the square-free part of $\Phi(E_\aa)$.

If $\aa$ and $\bb$ are two monomials in $E$ such that
$\supp(\aa)\subseteq\supp(\bb)$, then 
$$D(L_\bb)\otimes E(-\bb)\rTo^{\pi\otimes(\cdot \bb\aa^{-1})}
D(L_\aa)\otimes E(-\aa),$$ 
where $\pi$ is the map induced to divided powers by the 
canonical projection $\pi: L_\bb\rTo L_\aa$, is a morphism of
chain complexes lifting the inclusion
$E_\bb\subseteq E_\aa\subseteq E$.

Given a complex ${\mathbf F}_\bullet$ in $\B$, we may apply $\Phi$
to obtain a double complex of free $E$-modules, and we
set $\Phi({\mathbf F}_\bullet)$ to be the total complex of
this double complex. Because of the way $\Phi$ is defined on each
object of $\B_0$, this functor preserves minimality. The spectral
sequences of the double complex shows that it also preserves 
acyclicity.\Box

As Aramova, Avramov, Herzog [2000] and R\"omer [2001]
observe, the existence of such a construction
shows that if an $S$-module $M$
has a linear free resolution over $S$ if and only if the
corresponding $E$ module has a linear free resolution over $E$.
Our version of the construction also ``explains'' these authors'
formula for Betti numbers:

\corollary{betti} The  following equality holds among  Poincar\'e 
series:
$$\sum_{i=0}^\infty\sum_{\aa\in\NN^n}\beta_{i,\aa}^E(_EM)t^i\uu^\aa=
\sum_{i=0}^\infty\sum_{\aa\in\NN^n}\beta_{i,\aa}^S(_SM)
{{t^i\uu^\aa}\over{\prod_{j\in\supp(\aa)}(1-tu_j)}}$$
where $\beta_{i,\aa}^E(_EM)$ denotes the dimension of the 
degree $\aa$ part of $\tor^E_i(M,K)$, and similarly for
$\beta_{i,\aa}^S(_SM)$.

\references 
\parindent=0pt
\frenchspacing
\item{} A.~Aramova, L.A.~Avramov, J.~Herzog: 
Resolutions of monomial ideals and cohomology over exterior algebras, 
{\sl Trans. Amer. Math. Soc.} {\bf 352} (2000), no. 2, 579--594.
\medskip

\item{} A.~Aramova, J.~Herzog, T.~Hibi: 
Squarefree lexsegment ideals, 
{\sl Math. Z.} {\bf 228}, (1998), 353--378.
\medskip

\item{} D.~Bayer, H.~Charalambous, S.~Popescu:
Extremal Betti Numbers and  Applications to Monomial Ideals,
{\sl J. Algebra} {\bf 221}, (1999), 497--512.
\medskip

\item{}  L.~J.~Billera, J.~S.~Provan:  Decompositions of simplicial 
complexes related to diameters of convex polyhedra, 
{\sl Math. Oper. Res.} {\bf 5}, (1980), 576--594.
\medskip

\item{} A.~Bj\"orner: On the homology of geometric lattices, 
{\sl   Algebra Univ.}  {\bf 14}, (1982), 107--128.
\medskip

\item{} A.~Bj\"orner: The homology and shellability of matroids and
  geometric lattices, Chapter 7 of {\it Matroid Applications}, ed.
  Neil White, 226--283, Encyclopedia Math. Appl., {\bf 40}, Cambridge
  Univ. Press, Cambridge, 1992.  
\medskip
  
\item{} A.~Bj\"orner, G.~Ziegler: 
Combinatorial stratification of complex arrangements'', 
{\sl J. Amer. Math. Soc.} {\bf 5}, (1992), no. 1, 105--149. 
\medskip

\item{}
B.~Bollob\'as: {\it Modern Graph Theory}, Graduate Texts
in Mathematics {\bf 184}, Springer, New York, 1998.
\medskip

\item{} W.~Bruns, J.~Herzog: {\it Cohen-Macaulay Rings}, Cambridge
  Studies in advanced mathematics, {\bf 39}, Cambridge University
  Press 1993.  
\medskip
  
\item{} W.~Bruns, J.~Herzog: On multigraded resolutions, 
  {\sl Math. Proc. Camb. Phil. Soc.}, {\bf 118}, (1995), 245--257.  
\medskip

\item{} D.~Cohen, A.~Suciu: Characteristic varieties of arrangements,
{\sl Math. Proc. Cambridge Philos. Soc.} {\bf 127} (1999), no. 1, 33--53.
\medskip

\item{} H.~Crapo: A higher invariant for matroids, {\sl J. of
    Combinatorial Theory} {\bf 2}, (1967), 406--417.  
\medskip

\item{} J.~Eagon, V.~Reiner: Resolutions of Stanley-Reisner rings and 
Alexander duality, {\sl J. Pure Appl. Algebra} {\bf 130}, (1998), no.
3, 265--275. 
\medskip

\item{} D.~Eisenbud:
{\sl Commutative Algebra with a View Toward Algebraic Geometry},
Springer, New York, 1995.
\medskip

\item{}  D.~Eisenbud, F.-O.~Schreyer: 
Sheaf Cohomology and Free Resolutions over Exterior Algebras, 
preprint {\tt math.AG/0005055}.
\medskip

\item{} H.~Esnault, V.~Schechtman, E.~Viehweg:
Cohomology of local systems on the complement of hyperplanes, 
{\sl Invent. Math.} {\bf 109} (1992), no. 3, 557--561.
\medskip

\item{} M.~Falk: Arrangements and cohomology, 
{\sl Ann. Comb.} {\bf 1} (1997), no. 2, 135--157. 
\medskip

\item{} D.~Grayson, M.~Stillman: {\it Macaulay2}, a
software system devoted to supporting research in algebraic geometry
and commutative algebra.  Contact the authors, or download from
{\tt ftp://ftp.math.uiuc.edu/Macaulay2}.
\medskip

\item{}  M.~Hochster: Rings of invariants of tori, Cohen-Macaulay
rings generated by monomials, and polytopes, {\sl Ann. of Math.}
{\bf 96}, (1972), 318--337.
\medskip

\item{} M.~Hochster: Cohen-Macaulay rings, combinatorics and simplicial
complexes, in Ring theory II,  McDonald B.R., Morris, R. A. (eds),
{\it Lecture Notes in Pure and Appl. Math.} {\bf 26}, M. Dekker 1977.
\medskip

\item{}  A.~Libgober, S.~Yuzvinsky: Cohomology of the Orlik-Solomon algebras 
and local systems, {\it Compositio Math.} {\bf 121} (2000), no. 3,
337--361.
\medskip

\item{} W.~Massey, L.~Traldi:
On a conjecture of K. Murasugi,
{\sl Pacific J. Math.} {\bf 124} (1986),
no. 1, 193--213.
\medskip

\item{} D.~Matei, A.~Suciu: Cohomology rings and nilpotent 
quotients of real and complex arrangements, 
{\it Singularities and Arrangements, Sapporo-Tokyo 1998}, 
Advanced Studies in Pure Mathematics {\bf 27} (2000), 185--215.
\medskip

\item{} J.~Milnor:
Isotopy of links, in  {\it Algebraic geometry and topology.
A symposium in honor of S. Lefschetz}, pp. 280--306.
Princeton University Press, Princeton, N. J., 1957.
\medskip

\item{} M.~Musta\c t\v a: Local Cohomology at Monomial Ideals, 
in ``Symbolic computation in algebra, analysis, and geometry 
(Berkeley, CA, 1998)'', 
{\sl J. Symbolic Comput.} {\bf 29} (2000), no. 4-5, 709--720.
 \medskip

\item{} P.~Orlik, L.~Solomon: Combinatorics and topology of
complements of hyperplanes, {\sl Invent. Math.} {\bf 56}, (1980), 167--189.
\medskip

\item{} P.~Orlik, H.~Terao: {\it Arrangements of hyperplanes}, 
Grundlehren der Mathematischen Wissenschaften {\bf 300}, 
Springer-Verlag, Berlin, 1992.  
\medskip

\item{} J.~S.~Provan: Decompositions, shellings, and diameters of 
simplicial complexes and convex polyhedra, Thesis Cornell Univ. 1977.
\medskip

\item{} T.~R\"omer: Generalized Alexander Duality and Applications, 
Preprint 1999. To appear in {\sl Osaka J. Math.} {\bf 38} (2001).
\medskip

\item{} V.~Schechtman, H.~Terao, A.~Varchenko: Local systems over
complements of hyperplanes and the Kac-Kazhdan conditions for singular
vectors,  {\sl J. Pure Appl. Algebra} {\bf 100}, (1995), 93--102.
\medskip

\item{} R.~Stanley: {\it Combinatorics and Commutative Algebra}, Second
edition, Progress in Math. {\bf 41}, Birkh\"auser, 1996.
\medskip

\item{} R.~Stanley: Cohen-Macaulay rings and constructible
polytopes, {\sl Bull. Amer. Math. Soc.} {\bf 81}, (1975), 133-135. 
\medskip

\item{} N.~Terai: Generalization of Eagon-Reiner theorem
and $h$-vectors of graded rings, preprint 1997.     
\medskip

\item{} K.~Yanagawa: Alexander duality for Stanley-Reisner rings and
   square-free ${\bf N}^n$-graded modules, {\sl J. Algebra} 
{\bf 225} (2000), no. 2, 630--645. 
\medskip

\item{} S.~Yuzvinsky: Cohomology of the Brieskorn-Orlik-Solomon algebras, 
{\sl Comm. Algebra} {\bf 23}, (1995), 5339--5354.
\medskip

\item{} G.~Ziegler: On the difference  between real and 
complex arrangements,  {\sl Math. Z.} {\bf 212} (1993), 
no. 1,  1--11. 
\medskip

\bigskip\bigskip 
\vbox{\noindent Author Addresses: 
\smallskip 
\noindent{David Eisenbud}\par 
\noindent{Department of Mathematics, University of California, Berkeley, 
Berkeley CA 94720}\par 
\noindent{de@msri.org} 
\smallskip 
\noindent{Sorin Popescu}\par 
\noindent{Department of Mathematics, SUNY at Stony Brook,  
Stony Brook, NY 11794, and}\par 
\noindent{Department of Mathematics, Columbia  University, 
New York, NY 10027}\par 
\noindent{sorin@math.sunysb.edu}\par 
\smallskip 
\noindent{Sergey Yuzvinsky}\par 
\noindent{Department of Mathematics, University of Oregon, 
Eugene, OR 97403}\par 
\noindent{yuz@math.uoregon.edu}\par 
} 
 
\bye